\documentclass[final,3p,times]{elsarticle} 

\usepackage{amssymb}
\usepackage{amsthm}

\usepackage{amsmath,amsfonts}
\usepackage{array}
\usepackage{enumerate}
\usepackage{cases}
\usepackage{color}
\usepackage{here}

\newcommand{\Dev}[2]{{\frac{\partial {#1}}{\partial {#2}}}}

\newcommand{\IN}{\textrm{in}~~}
\newcommand{\ON}{\textrm{on}~~}

\newcommand{\Norm}[2]{\left\|#1\right\|_{#2}}
\newcommand{\Inner}[3]{\left(#1,#2\right)_{#3}}

\newcommand{\Set}[2]{{\left\{{#1}~;~{#2}\right\}}}

\newcommand{\be}{\begin{eqnarray}}
\newcommand{\ee}{\end{eqnarray}}
\newcommand{\bgn}{\begin{eqnarray*}}
\newcommand{\fin}{\end{eqnarray*}}

	\theoremstyle{definition} 
	\newcommand{\theoremname}{Theorem}
	\newtheorem{Theorem}{\theoremname}[section]
	\newcommand{\lemmaname}{Lemma}
	\newtheorem{lemma}[Theorem]{\lemmaname}
	\newcommand{\definitionname}{Definition}
	
	\newcommand{\assumptionname}{Assumption}

\begin{document}

\begin{frontmatter}



\title{Constructive a priori error estimates for a full discrete approximation of periodic solutions for the heat equation}


\author[label1]{Takuma Kimura\corref{cor1}}
\ead{tkimura@cc.saga-u.ac.jp}
\cortext[cor1]{Corresponding author}
\author[label1]{Teruya Minamoto}
\ead{minamoto@is.saga-u.ac.jp}
\author[label2]{Mitsuhiro T. Nakao}
\ead{mtnakao@imi.kyushu-u.ac.jp}
\address[label1]{Department of Information Science, Saga University, Saga 840-8502, Japan}
\address[label2]{Faculty of Science and Engineering, Waseda University, Tokyo 169-8555, Japan}

\begin{abstract}
 We consider the constructive a priori error estimates for a full discrete numerical solution of the heat equation with time-periodic condition. Our numerical scheme is based on the finite element semidiscretization in space direction combining with an interpolation in time by using the fundamental matrix for the semidiscretized problem. We derive the optimal order $H^1$ and $L^2$ error estimates, which play an important role in the numerical verification method of exact solutions for the nonlinear parabolic equations.  Several numeriacl examples which confirm us the optimal rate of convergence are presented. 
\end{abstract}

\begin{keyword}
Parabolic problem; 
Periodic solutions; 
Finite element method; 
Constructive a priori error estimates

\MSC 35B10 \sep 35K05 \sep 65M15 \sep 65M60
\end{keyword}

\end{frontmatter}



\section{Introduction}

Many works have been done concerning the error estimates for the approximate solutions of linear parabolic initial boundary value problems. Particularly, in \cite{Hansbo 1991}, \cite{Bernardi 1982}, they treated the time-periodic problems of the heat equation. On the other hand, recently, there are many results on the numerical enclosing the closed orbits corresponding to the periodic solutions  by mainly using  spectral techniques, \cite{Zgliczynski 2010},\cite{Figueras 2017} etc., as part of the study in dynamical systems. In their works, the spectral properties for the simple operator restricted to the rectangular domains are effectively used. In the present paper, we consider the finite element approach instead the spectral method. Such a technique seems to be more complicated and the error estimates are not so easy compared with spectral method. But, there is no limit to the shape of the domain at all. The method we describe here basically extends the results of the previous paper \cite{Nakao SIAM} to the time-periodic problem of a heat equation.

In the followings, we use the time-dependent Sobolev spaces with associated norms of the form $L^p((0,t); X)$. For example, $u \in L^2((0,T); H^1_0(\Omega))$, then
$$
\Norm{u}{L^2H^1_0}^2
\equiv
\Norm{u}{L^2((0,T); H^1_0(\Omega))}^2 
:=\int_0^T  \Norm{u(t)}{H^1_0(\Omega)}^2dt,
$$
also use the notation such that $\Norm{u}{L^2L^2} \equiv \Norm{u}{L^2((0,T); L^2(\Omega))}$ for short and so on. For other notations and properties of function spaces, see e.g. \cite{Barbu 1998}, \cite{Zeidler 1990}.          


\section{Problem and basic properties} 
\label{Problem}

In this section, we introduce the time-periodic problem and give the basic properties of the solution.\\
We consider the following heat equation with time-periodic condition:
\begin{subequations}
\begin{numcases}
{}
\displaystyle\frac{\partial u}{\partial t}-\nu\Delta u=f(x,t)
& $\IN \Omega \times J$, 
\label{eq:PSD} \\
u(x,t)=0
& $\ON \partial\Omega \times J$, 
\label{eq:PSD BC} \\
u(x,0)=u(x,T)
& $\IN \Omega$, \label{eq:PSD IC}
\end{numcases}\label{prob1}
\end{subequations}
where $\nu$ is a positive constant, $J:=(0,T)\subset\mathbb{R}~(T<\infty)$ and  $\Omega\subset\mathbb{R}^d$~$(d=1,2,3)$ a convex polygonal  or polyhedral domains. Also we define $\displaystyle \Delta_t\equiv\Dev{}{t}-\nu\Delta$ and assume that ${{f \in L^2((0,T); L^2(\Omega))}} \equiv L^2(\Omega \times J)$. On the existence and uniqueness of solution for \eqref{prob1}, see e.g. \cite{Barbu 1998}, \cite{Zeidler 1990}.

Now, for any $v \in L^2(\Omega)$ and  $t >0 $, we define the evolution operator $E(t): L^2(\Omega) \rightarrow L^2((0,t); H^1_0(\Omega))$ as a solution  $\phi \in L^2((0,t); H^1_0(\Omega))$ of the following equation. Namely, $E(t)v \equiv \phi$ satisfies
\begin{subequations}
\begin{numcases}
{}
\displaystyle\frac{\partial \phi}{\partial s}-\nu\Delta \phi=0
& $\IN \Omega \times (0,t)$, 
\label{eq:PSD-hom-right} \\
\phi(x,s)=0
& $\ON \partial\Omega \times (0,t)$, 
\label{eq:PSD-hom-right BC} \\
\phi(x,0)=v(x)
& $\IN \Omega$. \label{eq:PSD-hom-right IC}
\end{numcases}\label{hom-right}
\end{subequations}
Next, consider the solution $\psi \in L^2((0,t); H^1_0(\Omega))$ satisfying the following parabolic problem with homogeneous initial condition
\begin{subequations}
\begin{numcases}
{}
\displaystyle\frac{\partial \psi}{\partial s}-\nu\Delta \psi=f(x,s)
& $\IN \Omega \times (0,t)$, 
\label{eq:PSD-hom-int} \\
\psi(x,s)=0
& $\ON \partial\Omega \times (0,t)$, 
\label{eq:PSD-hom-int BC} \\
\psi(x,0)=0
& $\IN \Omega$. \label{eq:PSD-hom-int IC}
\end{numcases}\label{hom-int}
\end{subequations}
Then note that by using the notation in semigroup theory, e.g., \cite{Pazy 1983}, we can rewrite \eqref{hom-int} as follows:
\[
\psi(t)= \int_0^tE(t-s)f(s)ds.
\]
Taking notice that,  using a solution $\phi$ of \eqref{hom-right} for an appropriately chosen initial function $v=u(0)$ and $\psi$ in \eqref{hom-int}, the solution $u$ of \eqref{prob1} can be represented as $u(t) \equiv u(\cdot, t) = \phi(t) +\psi(t)$. Namely, we have 
\begin{align}
u(t) &= E(t)u(0) + \int_0^tE(t-s)f(s)ds \label{u(t)}.
\end{align}
Now, by the well known arguments using spectral theory in \cite{Barbu 1998} or semigroup approches in \cite{Pazy 1983}, for the minimal eigenvalue ${\lambda}_1$ of $-\Delta$ on $\Omega$, it holds that for the spaces $X = L^2(\Omega)$ or $X= H_0^1(\Omega)$
\begin{align}
  \Norm{E(t)v}{X}
    &\leq e^{-\nu {\lambda}_1 t}\Norm{v}{X}, \label{eq:Estimation}   
\end{align}
where $\Norm{u}{H_0^1(\Omega)} \equiv \Norm{\nabla u}{L^2(\Omega)}$. Then, from  the periodic condition, we have by \eqref{u(t)}
\begin{align}
  u(0) = E(T)u(0) + \psi(T).
   \label{eq:periodic-condition}   
\end{align}
Hence, from the contraction property of $E(T)$ due to the estimates \eqref{eq:Estimation}, the invertibility of the operator $I-E(T)$ follows and the initial value $u(0)$ is determined by
\begin{align}
  u(0) = (I-E(T))^{-1} \psi(T).
   \label{eq:u(0)}   
\end{align}
Furthermore, by the fact that $\psi$ is a solution of \eqref{hom-int}, it is readily seen that, by \eqref{eq:Estimation} and \eqref{eq:u(0)} (cf. in the proof of Lemma 4.2 of \cite{Nakao SIAM}): 
\begin{align}
   \Norm{u(0)}{L^2(\Omega)} \leq (1- e^{-\nu {\lambda}_1 T})^{-1} \frac{C_p}{\sqrt{\nu}}\Norm{f}{L^2L^2},
   \label{eq:L^2-norm-u(0)}   
\end{align}
where $C_p$ is a Poincar\'{e} constant on $\Omega$.
Also, if we use the fact that 
$\Norm{\psi(T)}{L^2(\Omega)} \leq \sqrt{T}\Norm{\psi_t}{L^2L^2}$ and the estimates $\Norm{\psi_t}{L^2L^2} \leq \Norm{f}{L^2L^2}$(Lemma 4.2 in \cite{Nakao SIAM}), we have another estimates as follows:
\begin{align}
   \Norm{u(0)}{L^2(\Omega)} \leq (1- e^{-\nu {\lambda}_1 T})^{-1} \sqrt{T}\Norm{f}{L^2L^2}.
   \label{eq:L^2-norm-u(0)-2}   
\end{align}

By the similar arguments, from \eqref{eq:Estimation}, \eqref{eq:u(0)} and the following estimates (cf. in the proof of Lemma 4.1 of \cite{Nakao SIAM})
$$
\Norm{\nabla \psi(T)}{L^2(\Omega)} \leq \frac{1}{\sqrt{\nu}}\Norm{f}{L^2L^2},
$$
we have the bound for $\nabla u(0)$ as
\begin{align}
   \Norm{\nabla u(0)}{L^2(\Omega)} \leq (1- e^{-\nu {\lambda}_1 T})^{-1} \frac{1}{\sqrt{\nu}}\Norm{f}{L^2L^2}.
   \label{eq:L^2-norm-grad-u(0)-2}   
\end{align}

The following lemma can be similarly obtained.
\begin{lemma}\label{Lem:estimates on Q}
For the solution $u$ of \eqref{prob1}, it holds that 
\begin{align}
   \Norm{u_t}{L^2L^2} &\leq \Norm{f}{L^2L^2},   \label{eq:u_t-norm on Q} \\
  \Norm{u(T)}{L^2(\Omega)}^2 + \nu \Norm{\nabla u}{L^2L^2}^2  &\leq (\frac{C_p^2}{\nu}+ T(1- e^{-\nu {\lambda}_1 T})^{-2})\Norm{f}{L^2L^2}^2. \label{eq:u(T)-grad-u-norm on Q} 
\end{align}
\end{lemma}
Proof. As in the proof of Lemma 4.1 in \cite{Nakao SIAM} we have
\be \Norm{u_t}{L^2(\Omega)}^2 + \nu\frac{d}{dt}\Norm{\nabla u}{L^2(\Omega)^d}^2
    \leq \Norm{f}{L^2(\Omega)}^2. \label{common-inequality}
\ee
Integrating this on $J$, by taking notice of the periodic condition, yields \eqref{eq:u_t-norm on Q}.\\
Similarly, from the proof of Lemma 4.2 in \cite{Nakao SIAM} we get
\begin{align}
\frac{d}{dt}\Norm{u(t)}{L^2(\Omega)}^2 + \nu\Norm{\nabla u(t)}{L^2(\Omega)^d}^2
  &\leq \frac{C_p^2}{\nu}\Norm{f}{L^2(\Omega)}^2, \label{common-inequality-2}
\end{align}
which proves \eqref{eq:u(T)-grad-u-norm on Q} by combining with the estimates \eqref{eq:L^2-norm-u(0)-2}. $\Box$


\section{Semidiscrete approximation}

In the present section, we define the semidiscrete approximation by the finite element method and derive the constructive error estimates. These results play important and essential roles in the error estimates for a full-discretization of the problem \eqref{prob1}. 

Let $S_h \equiv S_h(\Omega) \subset H^1_0(\Omega)$ be a finite dimensional subspace in spatial direction with $\dim S_h =n$ and let $V_k^1 \equiv V_k^1(J)\subset V^1(J)\equiv H^1(J) \cap \{u~|~u(0)=u(T)\}$ be a piecewise linear Lagrange type finite element space in time direction with $\dim S^k= m$. Also define $V:=H^1\bigl(J;L^2(\Omega)\bigr)\cap L^2\bigl(J;H_0^1(\Omega)\bigr)\cap \{u~|~u(0)=u(T)\; \text{in} \;H_0^1(\Omega)\}$. 

Now, let $P_h^1:H_0^1(\Omega)\to S_h$ be an $H_0^1$-projection  satisfying
\begin{align}
\Inner{\nabla (u-P_h^1u)}{\nabla v_h}{L^2(\Omega)}=0 \quad \forall v_h\in S_h, \label{eq:PSD:Definition of Ph1 in H01}
\end{align}
with the following assumptions on the approximation property:
\begin{align}
  \Norm{u-P_h^1u}{H_0^1(\Omega)}
    &\leq C_\Omega(h)\Norm{\Delta u}{L^2(\Omega)} ~\forall u\in H_0^1(\Omega) \cap Y(\Omega), \label{eq:IPP:I-Ph1 H01 error estimate} \\
  \Norm{u-P_h^1u}{L^2(\Omega)}
    &\leq C_\Omega(h)\Norm{u-P_h^1u}{H_0^1(\Omega)} ~\forall u\in H_0^1(\Omega). \label{eq:IPP:I-Ph1 L2 error estimate}
\end{align}
Here，$Y(\Omega):=\Set{u\in L^2(\Omega)}{\Delta u\in L^2(\Omega)}$. 

Now, we define the semidiscrete projection $P_h:V \to H^1\bigl(J;S_h(\Omega)\bigr) \equiv V^1\bigl(J;S_h(\Omega)\bigr)\cap  \{v_h(0)=v_h(T)\}$ by the following weak form: 
\begin{subequations}
\begin{numcases}
{}
\Inner{\tfrac{\partial}{\partial t}(u-P_hu)}{v_h}{L^2(\Omega)}+\nu\Inner{\nabla(u-P_hu)}{\nabla v_h}{L^2(\Omega)^d} = 0 ~\forall v_h\in S_h, t\in J, &\label{projection-P_h-1}\\
(P_hu)(0) = (P_hu)(T).  & \label{projection-P_h-2}
\end{numcases}\label{projection-P_h}
\end{subequations}
Note that $P_hu$ implies the semidiscrete approximation of a solution  $u$ for \eqref{prob1} with given function $f \in L^2(J;L^2(\Omega))$. Therefore,  we denote $(P_hu)(t)$ by $u_h(t)$, i.e., $u_h \equiv P_hu$ in the below.


Next we consider the constructive error estimates for $P_hu$  defined by \eqref {projection-P_h}.\\
For any $ v_h \in S_h$ and $t >0 $, we define the semidiscrete evolutional operator  $E_h(t):S_h \rightarrow S_h$ by the solution $\phi_h \in H^1\bigl((0,t);S_h(\Omega)\bigr)$ of the following equation. Namely, $E_h(t)v_h \equiv \phi_h$ corresponds to a semidiscretization of the solution $E(t)v \equiv \phi$ defined by \eqref{hom-right}.
\begin{subequations}
\begin{numcases}
{}
\displaystyle\frac{\partial \phi_h}{\partial s}-\nu\Delta_h \phi_h=0
& $\IN \Omega \times (0,t)$, 
\label{eq:discrete-hom-right} \\
\phi_h(x,0)=v_h(x)
& $\IN \Omega$. \label{eq:discrete-hom-right IC}
\end{numcases}\label{discrete-hom-right}
\end{subequations}
Here, $\Delta_h$ means the discretization of a weak Laplacian on $S_h$ and  \eqref{eq:discrete-hom-right} is equivalent to the following variational form:
\begin{eqnarray}
((\phi_h)_t, \eta_h)_{L^2(\Omega)} + \nu (\nabla \phi_h, \nabla \eta_h)_{L^2(\Omega)} = 0 \quad \forall \eta_h \in S_h, \; t>0.
\label{eq:discrete-var.}
\end{eqnarray}
Similarly, as an semidiscretization for \eqref{hom-int}, we consider a solution $\psi_h \in H^1\bigl((0,t);S_h(\Omega)\bigr)$ of the following  
equation
\begin{subequations}
\begin{numcases}
{}
\displaystyle\frac{\partial \psi_h}{\partial s}-\nu\Delta_h \psi_h=P_h^0f
& $\IN \Omega \times (0,t)$, 
\label{eq:discrete-hom-int} \\
\psi_h(x,0)=0
& $\IN \Omega$, \label{eq:discrete-hom-int IC}
\end{numcases}\label{discrete-hom-int}
\end{subequations}
where $P_h^0f$ means the $L^2$-projection of $f$ to $S_h$.
 Also by using the similar symbol and arguments as in the previous section we get the following expression:
\begin{eqnarray}
P_hu(t) = E_h(t)u_h(0) + \int_0^tE_h(t-s)P_h^0fds.
\label{eq:semi-discrete}
\end{eqnarray}
Here, note that we can numerically  compute the norm 
 $\kappa_1 :=\Norm{E_h(T)}{\mathcal{L}(H^1_0)}$ by matrix norm computations to confirm it is actually less than one, namely, contraction map on $S_h$. On the actual estimation of $\kappa_1$, see Remark 4.1 in the next section. And we can also compute the following inverse operator norm for $(I-E_h(T))^{-1}$
\be
\Norm{(I-E_h(T))^{-1}}{\mathcal{L}(H^1_0)}  \leq  (1- \kappa_1)^{-1}. \label{kappa_1}
\ee
Thus, from the definition and discrete analog to the previous section, we have $u_h(0) =  (I-E_h(T))^{-1}\psi_h(T)$ and obtain the following estimates:
\be
\Norm{\nabla u_h(0)}{L^2(\Omega)} \leq (1- \kappa_1)^{-1}\frac{1}{\sqrt{\nu}}\Norm{f}{L^2L^2}. \label{est-grad-u_h(0)}
\ee

Now, in order to get the error estimates for the semidisctrete approximation defined by \eqref{projection-P_h} or equivalently  by \eqref{eq:semi-discrete} for the problem  \eqref{prob1}, first we consider the constructive error estimates for the semidiscretization of the nonhomogeneous parabolic initial boundary value problem with initial condition $\xi_0 \in H^1_0(\Omega)$ of the form :
\begin{subequations}
\begin{numcases}
{}
\displaystyle\frac{\partial \xi}{\partial t}-\nu\Delta \xi=f(x,t)
& $\IN \Omega \times J$, 
\label{nonhom} \\
\xi(x,t)=0
& $\ON \partial\Omega \times J$, 
\label{nonhom BC} \\
\xi(x,0)=\xi_0
& $\IN \Omega$. \label{nonhom IC}
\end{numcases}\label{nonhom-parab}
\end{subequations}
Let $\xi_h \in S_h$ be a semidiscrete approximation of \eqref{nonhom-parab} given by the following weak form:
\begin{subequations}
\begin{numcases}
{}
((\xi_h)_t, v_h)_{L^2(\Omega)} + \nu (\nabla \xi_h, \nabla v_h)_{L^2(\Omega)} = (f(\cdot, t),v_h)_{L^2(\Omega)}  \; \forall v_h \in S_h, \; t>0& \label{nonhom-semidiscrete 1}\\
\xi_h(0) = \zeta_h. &
\label{nonhom-semidiscrete 2}
\end{numcases}\label{nonhom-semidiscrete}
\end{subequations}
Here, $\zeta_h \in S_h$ is an appropriate approximation of $\xi_0$. Then we have the following estimates for solutions of \eqref{nonhom-parab} and \eqref{nonhom-semidiscrete}.
\begin{lemma}\label{est-u-u_h-Q}
\begin{align}
\Norm{\xi_t}{L^2L^2} &\leq \Norm{f}{L^2L^2}+\sqrt{\nu}\Norm{\nabla \xi_0}{L^2(\Omega)} \label{est-u_t-Q}, \\
\Norm{\xi}{L^2H^1_0} &\leq \frac{C_p}{\nu}\Norm{f}{L^2L^2}+\frac{1}{\sqrt{\nu}}\Norm{\xi_0}{L^2(\Omega)} \label{est-grad-u-Q}, \\
\Norm{(\xi_h)_t}{L^2L^2} &\leq \Norm{f}{L^2L^2}+\sqrt{\nu}\Norm{\nabla \zeta_h}{L^2(\Omega)} \label{est-u_h-t-Q}, \\
\Norm{\xi_h}{L^2H^1_0} &\leq \frac{C_p}{\nu}\Norm{f}{L^2L^2}+\frac{1}{\sqrt{\nu}}\Norm{\zeta_h}{L^2(\Omega)}. \label{est-grad-u_h-Q}
\end{align}
\end{lemma}
Proof.　These results are obtained by the similar arguments to that in the proofs for Lemma 4.1-4.4 in \cite{Nakao SIAM} with some additional considerations. \\First, by the same argument to derive \eqref{common-inequality}, we have 
\be \Norm{\xi_t}{L^2L^2}^2 + \nu\Norm{\nabla \xi(T)}{L^2(\Omega)}^2
    \leq \Norm{f}{L^2L^2}^2 +\nu\Norm{\nabla{\xi}(0)}{L^2(\Omega)}^2, \label{common-inequality-1}
\ee
which implies  \eqref{est-u_t-Q}. 
Next, by the similar manner of getting \eqref{common-inequality-2} in the proof of Lemma \ref{Lem:estimates on Q}, we have
\begin{align*}
 \frac{d}{dt}\Norm{\xi(t)}{L^2(\Omega)}^2 + \nu\Norm{\nabla \xi(t)}{L^2(\Omega)^d}^2
     &\leq \frac{C_p^2}{\nu}\Norm{f}{L^2(\Omega)}^2.
\end{align*}
Thus integrating both sides in $t$ yields the estimates \eqref{est-grad-u-Q}.

We now take $v_h := (\xi_h)_t$ for $t>0$ in \eqref{nonhom-semidiscrete 1} and integrate it in $t$, we have
\begin{align}
  \Norm{(\xi_h)_t}{L^2L^2}^2 + \nu\Norm{\nabla \xi_h(T)}{L^2(\Omega)^d}^2
    &\leq \Norm{f}{L^2L^2}^2 + \nu\Norm{\nabla \xi_h(0)}{L^2(\Omega)^d}^2, \label{common-inequality-3}
\end{align}  
which proves the assertion \eqref{est-u_h-t-Q}. Finally, the estimates \eqref{est-grad-u_h-Q} can be easily derived by the argument analogous to proving \eqref{est-grad-u-Q}. $\Box$

Also, setting $\xi_{\perp} := \xi- \xi_h$, we obtain the following two kinds of  error estimates, which are obtained similar arguments in the proof of Theorem 4.6 in \cite{Nakao SIAM}. 
\begin{Theorem}\label{H^1_0-error-est-non-hom}
The following estimates for $\xi_{\perp} := \xi- \xi_h$ hold: 
\begin{align}
  \Norm{\xi_{\perp}}{L^2H^1_0} \leq \bigl[ \frac{C_\Omega(h)^2}{\nu^2} \bigl\{4\Norm{f}{L^2L^2}^2 + \nu(\Norm{\nabla \xi_0}{L^2(\Omega)}^2+\Norm{\nabla \zeta_h}{L^2(\Omega)}^2)\bigr\} + \frac{1}{2\nu}\Norm{\xi_0 - \zeta_h}{L^2(\Omega)}^2\bigr]^{\frac{1}{2}}, \label{eq:H^1_0-error-est-non-hom}
\end{align}
also $L^2$-estimates at $t=T$,
\begin{align}
  \Norm{\xi_{\perp}(T)}{L^2(\Omega)} \leq \bigl[\frac{2}{\nu}C_\Omega(h)^2 \bigl\{4\Norm{f}{L^2L^2}^2 + \nu(\Norm{\nabla \xi_0}{L^2(\Omega)}^2+\Norm{\nabla \zeta_h}{L^2(\Omega)}^2)\bigr\} + \Norm{\xi_0 - \zeta_h}{L^2(\Omega)}^2\bigr]^{\frac{1}{2}}.  \label{eq:L^2-error-est-non-hom-t=T}
\end{align}
\end{Theorem}

Proof. Applying the same arguments in the proof of Theorem 4.6 in \cite{Nakao SIAM}, we have 
\begin{align*}
  \frac{1}{2}\frac{d}{dt}\Norm{\xi_\perp}{L^2(\Omega)}^2 + \nu\Norm{\xi_\perp}{H_0^1(\Omega)}^2
       &\leq \frac{C_\Omega(h)^2}{\nu}\left(2\Norm{f}{L^2(\Omega)}^2 + \Norm{\frac{\partial \xi_h}{\partial t}}{L^2(\Omega)}^2 + \Norm{\frac{\partial \xi}{\partial t}}{L^2(\Omega)}^2\right).
\end{align*}
Integrating this on $J$, from \eqref{common-inequality-1} and \eqref{common-inequality-3}, we get
\begin{align*}
\frac{1}{2}\Norm{\xi_\perp(T)}{L^2(\Omega)}^2 + \nu\Norm{\xi_\perp}{L^2H_0^1}^2
  &\leq \frac{C_\Omega(h)^2}{\nu}\left(2\Norm{f}{L^2L^2}^2 + \Norm{\frac{\partial \xi_h}{\partial t}}{L^2L^2}^2 + \Norm{\frac{\partial \xi}{\partial t}}{L^2L^2}^2\right) + \frac{1}{2}\Norm{\xi_\perp(0)}{L^2(\Omega)}^2 \\
    &\leq\frac{C_\Omega(h)^2}{\nu}\left(4\Norm{f}{L^2\bigl(J;L^2(\Omega)\bigr)}^2 + \nu\Norm{\nabla \xi_h(0)}{L^2(\Omega)^d}^2 + \nu\Norm{\nabla \xi(0)}{L^2(\Omega)^d}^2 \right) + \frac{1}{2}\Norm{\xi_0-\zeta_h}{L^2(\Omega)}^2,
\end{align*}
which yields the desired conclusions \eqref{eq:H^1_0-error-est-non-hom} and \eqref{eq:L^2-error-est-non-hom-t=T}. $\Box$

\section{Full-discrete approximation and error estimates}
In this section, we define the full-discrete approximation of solutions for the problem \eqref{prob1} by using an interpolation procedure in time direction for the spatial discretized solution. We also show a computational scheme for this full discretization by the effective use of the fundamental matrix for an ODE system corresponding to semidiscretized problem. The constructive and optimal order $H^1$ and $L^2$ error estimates are established, which are main results in the present paper.

\subsection{A full discretizaion scheme}

Now, defining the interpolation operator $\Pi^k:$ $~V^1(J)$ $\to$ $V_k^1$ in time direction by
$$u(t_i)=\Pi^ku(t_i) \quad ^\forall i\in\{0,1,\cdots,m\},$$
we define the full discrete projection $P_h^k:V\to V_k^1\bigl(J;S_h(\Omega)\bigr) \equiv S_h\otimes V_k^1$ as
\begin{eqnarray}
P_h^ku:=\Pi^k(P_hu),
\label{full-discrete-scheme}
\end{eqnarray}
which corresponds to the full discretization of \eqref{prob1}.

In order to present the actual computation procedure of the above full discretization scheme,  we first consider a representation of the semidiscretization defined in \eqref{projection-P_h}.    Let  $\{\phi_i\}_{i=1}^n$ be a basis of $S_h$ and define the $n \times n$ matrices $L_\phi$, $D_\phi$ by 
\begin{align}
L_{\phi,i,j}&:=\Inner{\phi_j}{\phi_i}{L^2(\Omega)}, \qquad D_{\phi,i,j}:=\Inner{\nabla \phi_j}{\nabla \phi_i}{L^2(\Omega)^d}, \label{eq:PSD:Definition of Lphi and Dphi} 
\end{align}
respectively. Since they are symmetric and positive definite, we get the  Cholesky decomposition as $L_\phi=L_\phi^{1/2}L_\phi^{T/2}$ and $D_\phi=D_\phi^{1/2}D_\phi^{T/2}$, respectively. 
Also note that there exists a vector valued function $\vec{u}_h\in V^1(J)^n$ satisfying 
$$P_hu(x,t)=\vec{u}_h(t)^T\Phi(x),$$
where $\Phi(x) \equiv (\phi_1,\cdots, \phi_n)^T$. \\
Thus by using 
$\vec{u}_h$，the semidiscretization \eqref{projection-P_h} is equivalently presented as ODEs:
\begin{subequations}
\begin{numcases}
{}
{L_\phi}\frac{d}{dt}\vec{u}_h+\nu {D_\phi}\vec{u}_h=\tilde{f}
& $\IN J$, 
\label{eq:temp2 eq} \\
\vec{u}_h(0)=\vec{u}_h(T),
&
\label{eq:temp2 IC}
\end{numcases}\label{eq:temp2}
\end{subequations}
where $\tilde{f}  = (\tilde{f}_i) \in\mathbb{R}^n$ with $\tilde{f}_i=\Inner{f}{\phi_i}{L^2(\Omega)}.
$
For simplicity we denote as ${\vec{b}}(t) \equiv L_\phi^{-1}\tilde{f}(t)$.  Then note that using the fundamental matrix $\Theta(t)=\exp(-\nu L_\phi^{-1}D_\phi t)$ of the equation \eqref{eq:temp2 eq}, we can represent 
\eqref{eq:temp2} as
\begin{subequations}
\begin{numcases}
{}
{\vec{u}_h}(t)=\Theta(t)\vec{u}_h(0)+\int_0^t\Theta(t-s){\vec{b}}(s)~ds
& $\IN J$, 
\label{eq:temp3 eq} \\
\vec{u}_h(0)=\vec{u}_h(T).
&
\label{eq:temp3 IC}
\end{numcases}\label{eq:temp3}
\end{subequations}
Therefore, assuming that the invertibility of  $\left(I-\Theta(T)\right)$，from \eqref{eq:temp3 eq}，we have
\begin{eqnarray*}
\vec{u}_h(0)=\vec{u}_h(T)
&\Longleftrightarrow&
{\vec{u}_h}(0)=\Theta(T)\vec{u}_h(0)+\int_0^T\Theta(T-s){\vec{b}}(s)~ds,
\\
&\Longleftrightarrow&
{\vec{u}_h}(0)=\left(I-\Theta(T)\right)^{-1}\int_0^T\Theta(T-s){\vec{b}}(s)~ds,
\end{eqnarray*}
which yields the following expression of the solution of \eqref{eq:temp3}：
\begin{eqnarray}
{\vec{u}_h}(t)=\Theta(t)\left(I-\Theta(T)\right)^{-1}\int_0^T\Theta(T-s){\vec{b}}(s)~ds+\int_0^t\Theta(t-s){\vec{b}}(s)~ds.
\label{eq:Sol.ODEs}
\end{eqnarray}
Hence, we obtain 
$$
P_h^ku(x,t_j)=\left(\Theta(t_j)\left(I-\Theta(T)\right)^{-1}\int_0^T\Theta(T-s){\vec{b}}(s)~ds+\int_0^{t_j}\Theta(t_j-s){\vec{b}}(s)~ds\right) \cdot \Phi(x).
$$
Thus the full discrete approximation $P_h^k u \equiv \Pi^kP_h u$ for the solution  $u$ of \eqref{prob1} can be numerically computed by using this procedure. \\ 

\noindent
{\it Remark 4.1}:

For any $v_h \in S_h$, using the definition \eqref{eq:PSD:Definition of Lphi and Dphi}, by some simple consideration on the $H^1_0$ norm for the element  $E_h(T)v_h \in S_h$,  we have readily seen that 
 \begin{align*}
  \Norm{E_h(T)v_h}{H^1_0}
    & \leq \Norm{D_\phi^{T/2}\exp(-\nu T L_\phi^{-1}D_\phi)D_\phi^{-T/2}}{2}\Norm{v_h}{H^1_0}, 
   \end{align*}
where $||\cdot||_2$ means the matrix 2-norm. This immediately yields the estimate of $\kappa_1$ in \eqref{kappa_1}.

\subsection{$H^1$ error estimates}
\label{H^1-error}

In this subsecton, we present an error estimate in the $L^2H^1_0$ sense on $ \Omega  \times J$ for the full discretization \eqref{full-discrete-scheme}. Denoting again the semidiscrete projection $P_hu$ defined in \eqref{projection-P_h} as $P_hu \equiv u_h$, the semidiscrete approximation  $u_h$ for \eqref{prob1} is written by
\begin{subequations}
\begin{numcases}
{}
\displaystyle\frac{\partial u_h}{\partial t}-\nu\Delta_h u_h=P_h^0f
& $\IN \Omega \times J$, 
\label{semidiscrete-int} \\
u_h(\cdot,0)=u_h(\cdot,T)
& $\IN \Omega$. \label{semidiscrete-int IC}
\end{numcases}\label{semidiscrete}
\end{subequations}
In order to obtain the desired estimates, we use the following  decomposition
\be
u-P_h^ku = (u-u_h) + (u_h -\Pi^ku_h). \label{full-error-est}
\ee
The second term of the above is estimated by using the standard interpolation estimates, e.g., \cite{Schultz 1973}, we have from \eqref{est-u_h-t-Q} and \eqref{est-grad-u_h(0)} 
\begin{align}
  \Norm{u_h-\Pi^ku_h}{L^2L^2}
    &\leq C_J(k) \Norm{(u_h)_t}{L^2L^2} \nonumber\\
    &\leq  C_J(k) (\Norm{f}{L^2L^2} +\sqrt{\nu}\Norm{\nabla (u_h(0))}{L^2(\Omega)}) \nonumber\\
    &\leq  C_J(k) (\Norm{f}{L^2L^2} +(1-\kappa_1)^{-1}\Norm{f}{{L^2L^2}}) \nonumber\\
    &=  C_J(k) \frac{2-\kappa_1}{1-\kappa_1}\Norm{f}{{L^2L^2}}. \label{u_h-P_h^ku}
\end{align}
Furthermore, using an inverse estimation constant $C_{inv}(h)$, which makes possible to bound the $H^1$ norm by the $L^2$ norm in $S_h$, we get 
\be
\Norm{u_h-\Pi^ku_h}{L^2H^1_0} \leq C_{inv}(h)C_J(k) \frac{2-\kappa_1}{1-\kappa_1}\Norm{f}{L^2L^2}. \label{u_h-P_h^ku_H^1}
\ee
Note that using the definition of the operator $E_h(t)$, we have by \eqref{semidiscrete}
\be
u_h(0) = E_h(T)u_h(0) + \psi_h(T).
\ee
{{Therefore, using $\psi(t)$ defined by \eqref {hom-int},}} we have
\begin{align}
  u(0)-u_h(0)
    & = E(T)u(0) +  \psi(T) -(  E_h(T)u_h(0) + \psi_h(T)) \nonumber\\
    & = E(T)(u(0) -u_h(0)) + (E(T)-  E_h(T))u_h(0) + (\psi(T)- \psi_h(T)), \nonumber
\end{align}
which implies
\begin{align}
  (I-E(T))(u(0)-u_h(0))
    & =(E(T)u_h(0) + \psi(T)) -  (E_h(T)u_h(0) + \psi_h(T)). \label{right-side}
   \end{align}
Note that, for any $t \in J$,  setting 
$$
\xi(t) := E(t)u_h(0) + \psi(t), \; \xi_0 :=u_h(0)
$$
$$
\xi_h(t) := E_h(t)u_h(0) + \psi_h(t), \; \zeta_h :=u_h(0),
$$
then $\xi$ and $\xi_h$ are solutions corresponding to \eqref{nonhom-parab} and \eqref{nonhom-semidiscrete}, respectively. \\
Hence, setting $\xi_{\perp} :=\xi- \xi_h $, the right-hand side of \eqref{right-side} coincides with $\xi_{\perp}(T)$. 
Therefore, we have
\begin{eqnarray}
  \Norm{u(0)-u_h(0)}{L^2(\Omega)}
  &=&\Norm{(I-E(T))^{-1}\xi_{\perp}(T)}{L^2(\Omega)}\notag\\
  &\leq& \Norm{(I-E(T))^{-1}}{\mathcal{L}(L^2)} \Norm{\xi_{\perp}(T)}{L^2(\Omega)}. \label{u(0)-u_h(0)-est}
   \end{eqnarray}
By the argument in the section \ref{Problem}, we have the following estimates
\begin{align}
  \Norm{(I-E(T))^{-1}}{\mathcal{L}(L^2)}
    & \leq (1- e^{-\nu {\lambda}_1 T})^{-1}. \label{first-est}
   \end{align}
Next, applying the error estimates \eqref{eq:L^2-error-est-non-hom-t=T} in Theorem \ref{H^1_0-error-est-non-hom} with taking notice of $\xi_0 = \zeta_h$, by using \eqref{est-grad-u_h(0)} we have
\begin{align}
  \Norm{\xi_{\perp}(T)}{L^2(\Omega)}
    & \leq \{ \frac{2}{\nu}C_\Omega(h)^2 \bigl(4\Norm{f}{L^2L^2}^2 + \nu(2\times \Norm{\nabla u_h(0)}{L^2}^2 \bigr) + 0\}^{\frac{1}{2}} \nonumber \\
    & = \frac{2}{\sqrt{\nu}}C_\Omega(h) \bigl(2+(1-\kappa_1)^{-2})^{\frac{1}{2}}\Norm{f}{L^2L^2}.   \label{second-est}
   \end{align}
Therefore, from \eqref{u(0)-u_h(0)-est}-\eqref{second-est}, we obtain
\begin{align}
  \Norm{u(0)-u_h(0)}{L^2(\Omega)}
    & \leq K_1C_\Omega(h) \Norm{f}{L^2L^2},   \label{u(0)-u_h(0)-norm}
   \end{align}
where
$$
K_1 \equiv \frac{2}{\sqrt{\nu}}(1- e^{-\nu {\lambda}_1 T})^{-1} \bigl(2+(1-\kappa_1)^{-2})^{\frac{1}{2}}.
$$
On the other hand, we have by \eqref{eq:H^1_0-error-est-non-hom} in Theorem \ref{H^1_0-error-est-non-hom}
\be
\Norm{u-u_h}{L^2H^1_0} \leq \bigl\{ \frac{C_\Omega(h)^2}{\nu^2} \bigl(4\Norm{f}{L^2L^2}^2&  +& \nu(\Norm{\nabla u(0)}{L^2(\Omega)}^2 
+\Norm{\nabla u_h(0)}{L^2(\Omega)}^2) \nonumber \\
 & + &\frac{1}{2\nu}\Norm{u(0) - u_h(0)}{L^2(\Omega)}^2\bigr)\bigr\}^{\frac{1}{2}}. \label{u-u_h}
\ee 
Thus, from the estimates \eqref{eq:L^2-norm-grad-u(0)-2}, \eqref{est-grad-u_h(0)} and \eqref{u(0)-u_h(0)-norm}, we obtain the following estimation for the semidiscrete solution:
\be
\Norm{u-u_h}{L^2H^1_0} \leq K_2C_\Omega(h)\Norm{f}{L^2L^2}, \label{error-est-semidiscrete}
\ee
where we set as
$$ 
K_2 \equiv \frac{1}{\nu}\bigl\{4+5(1- e^{-\nu {\lambda}_1 T})^{-2}+(1+(1- e^{-\nu {\lambda}_1 T})^{-2})(1-\kappa_1)^{-2}\bigr\}^{\frac{1}{2}}.
$$
Combining \eqref{u_h-P_h^ku_H^1} and \eqref{error-est-semidiscrete} with \eqref{full-error-est}, we have the following desired $H^1$ error estimates.
\begin{Theorem}\label{H^1_0-error-est}
Let $P_h^ku$ be a full-discrete approximation defined by \eqref{full-discrete-scheme} for the  periodic solution $u$ of the heat equation \eqref{prob1}. Then, it holds that
\be
\Norm{u-P_h^ku}{L^2H^1_0} \leq \{K_2C_\Omega(h)+ C_{inv}(h)C_J(k) \frac{2-\kappa_1}{1-\kappa_1}\}\Norm{f}{L^2L^2}. 
\ee
\end{Theorem}
Here, the constant $K_2$ is defined in \eqref{error-est-semidiscrete}.

\subsection{$L^2$ error estimates}

In this subsection, we consider the error estimates in the $L^2L^2$ sense for the full-discrete approximation $P_h^ku$, which enable us higher order estimates than  the $L^2H^1$ error bound in Theorem \ref{H^1_0-error-est}. As in the previous subsection, we use a semidiscrete approximation $u_h$ with decomposition \eqref{full-error-est}. Note that, if we take as  $k \approx h^2$, then by applying the $L^2$ estimates \eqref{u_h-P_h^ku}, we immediately obtain $O(h^2)$ estimates for the latter term in \eqref{full-error-est}. Hence,  it suffices to derive the  $O(h^2)$ estimates for the former part.

\begin{Theorem}\label{L2-error-est}
It holds that  
\begin{align}
  \Norm{u-P_h^ku}{L^2L^2} \leq\{(\frac{3-2\kappa_1}{1-\kappa_1}\frac{2}{\nu}+2K_2)C_\Omega(h)^2 +\frac{2-\kappa_1}{1-\kappa_1} C_J(k) \}\Norm{f}{L^2L^2},
\end{align}
where $K_2$ is the same constant defined  in the estimates \eqref{error-est-semidiscrete}.
\end{Theorem}
Proof.　For any function $g \in L^2(Q)$, where $Q \equiv \Omega \times J$, let $v$ be a solution of \eqref{prob1} with the right-hand side $g(T-t) \equiv g(\cdot,T-t)$. Here, $t$ is a variable such that $t \in J$. Then $v$ satisfies the following weak form:
\be
\Inner{\tfrac{\partial}{\partial t}v(t)}{w}{L^2(\Omega)}+\nu\Inner{\nabla v(t)}{\nabla w}{L^2(\Omega)^d} = \Inner{g(T-t)}{w}{L^2(\Omega)} & \; \forall w \in H^1_0(\Omega), t\in J. \label{backward-weak-form}
\ee
Particularly, taking $w=u-u_h$ in \eqref{backward-weak-form} and transform the variable as $t \rightarrow T-s$, we have 
\bgn
\Inner{-\tfrac{\partial}{\partial s}v(T-s)}{u-u_h}{L^2(\Omega)}+\nu\Inner{\nabla (u-u_h)}{\nabla v(T-s)}{L^2(\Omega)^d} = \Inner{g(s)}{ u-u_h}{L^2(\Omega)}. \label{weak-form-1}
\fin
Integrating both sides of the above in $s$ on $(0,T)$ yields that
\be
\int_0^T\Inner{-\tfrac{\partial}{\partial s}v(T-s)}{u-u_h}{L^2(\Omega)}ds+\nu\int_0^T\Inner{\nabla (u-u_h)}{\nabla v(T-s)}{L^2(\Omega)^d}ds & \nonumber \\ 
 = \Inner{g(s)}{u-u_h}{L^2(Q)}. \label{weak-form-2}
\ee
Taking notice of the periodic condition, observe that
$$
\int_0^T\Inner{-\tfrac{\partial}{\partial s}v(T-s)}{u-u_h}{L^2(\Omega)}ds=
\int_0^T\Inner{\tfrac{\partial}{\partial s}(u-u_h)}{v(T-s)}{L^2(\Omega)}ds.
$$
Therefore, by the definition of $u_h$ and \eqref{weak-form-2} we have for any $v_h(s) \in S_h$
\begin{eqnarray}
  \Inner{g}{u-u_h}{L^2(Q)} 
& = & 
\int_0^T\Inner{\tfrac{\partial}{\partial s}(u - u_h)}{(v-v_h)(T-s)}{L^2(\Omega)}ds 
\nonumber\\
&&+ \nu\int_0^T\Inner{\nabla (u-u_h)}{\nabla (v-v_h)(T-s)}{L^2(\Omega)^d}ds
\nonumber\\
&\leq &\Norm{(u-u_h)_t}{L^2L^2}\Norm{v-v_h}{L^2L^2} 
\nonumber \\
&&+\nu\Norm{\nabla(u-u_h)}{L^2L^2}\Norm{\nabla(v-v_h)}{L^2L^2}.
\label{evaluation-1}
\end{eqnarray}
Moreover, by the similar derivation process of \eqref{u_h-P_h^ku} using \eqref{est-u_h-t-Q} in the previous subsection and Lemma \ref{Lem:estimates on Q}, we obtain
\be
\Norm{(u-u_h)_t}{L^2L^2} & \leq & \Norm{u_t}{L^2L^2} + \Norm{(u_h)_t}{L^2L^2} \nonumber \\ 
& \leq & \Norm{f}{L^2L^2}+\frac{2-\kappa_1}{1-\kappa_1}\Norm{f}{{L^2L^2}} \nonumber \\
& = & \frac{3-2\kappa_1}{1-\kappa_1}\Norm{f}{{L^2L^2}}. \label{L2-1}
\ee
Furthermore, for any $t \in (0,T)$, taking $v_h(t):=P_h^1v(t)$ to apply the approximation properties  \eqref{eq:IPP:I-Ph1 H01 error estimate} and  \eqref{eq:IPP:I-Ph1 L2 error estimate}, by considering the estimates in Lemma \ref{Lem:estimates on Q}, we have
\be
\Norm{\nabla(v-v_h)}{L^2L^2} & \leq & C_\Omega(h)\Norm{\Delta v}{L^2L^2} \nonumber \\
& \leq & C_\Omega(h)\frac{1}{\nu}\Norm{v_t -g}{L^2L^2} \nonumber \\
& \leq & C_\Omega(h)\frac{1}{\nu}(\Norm{v_t}{L^2L^2}+\Norm{g}{L^2L^2})\nonumber \\
& \leq & \frac{2}{\nu} C_\Omega(h)\Norm{g}{L^2L^2} \label{L2-3}
\ee
and
\be
\Norm{v-v_h}{L^2L^2} \leq \frac{2}{\nu} C_\Omega(h)^2\Norm{g}{L^2L^2}. \label{L2-4}
\ee
Therefore, combining \eqref{L2-1}-\eqref{L2-4} with \eqref{error-est-semidiscrete},
 we have the estimates
\be
\Norm{u-u_h}{L^2L^2} \leq (\frac{3-2\kappa_1}{1-\kappa_1}\frac{2}{\nu}+2K_2)C_\Omega(h)^2 \Norm{f}{L^2L^2},
\ee
which proves the theorem by \eqref{u_h-P_h^ku}. $\Box$
%
\section{Numerical examples}
In this section, we show several numerical examples which confirm us the optimal rate of convergence. 
We used the interval arithmetic toolbox INTLAB 11~\cite{INTLAB} with MATLAB R2012a on an Intel Xeon W2155 (3.30 GHz) with CentOS 7.4. 

Here, we only consider $d=1$, $\Omega=(0,1)$ and $J=(0,1)$, 
then the lower bound of eigenvalue of $-\Delta$ on $\Omega$ can be taken as $\lambda_1=\pi^2$. 
Furthermore, we set $f$ to be the problem \eqref{prob1} have the exact solution $u(x,t)=\sin(2\pi x)\sin(2\pi t+\beta)$. 
Here, $\beta$ is a given constant. 
Since the exact solutions are known, the upper bounds of the exact errors for approximate solutions can be validated in the a posteriori sense. 

We used the finite dimensional subspaces $S_h$ and $V_k^1$ spanned by piecewise
linear basis functions with uniform mesh size $h$ and $k$, respectively. 
Therefore, the constants can be taken as 
$C_\Omega(h) = h/\pi$, $C_J(k)=k/\pi$, ${{C_{inv}(h)}} =\sqrt{12}/h$, and $C_p=1/\pi$, respectively. 
We set $k=h^2$ then Theorem 4.1 and 4.2 are $O(h)$ and $O(h^2)$ error estimates . 
In Figure 1-2, {{the a priori error estimates and the}} exact errors of this example are shown. 
These Figures show the estimates presented in Theorem 4.1-4.2 give the optimal order estimates. 
\begin{figure}[htbp]
\begin{center}
\includegraphics[clip,width=0.95\linewidth]{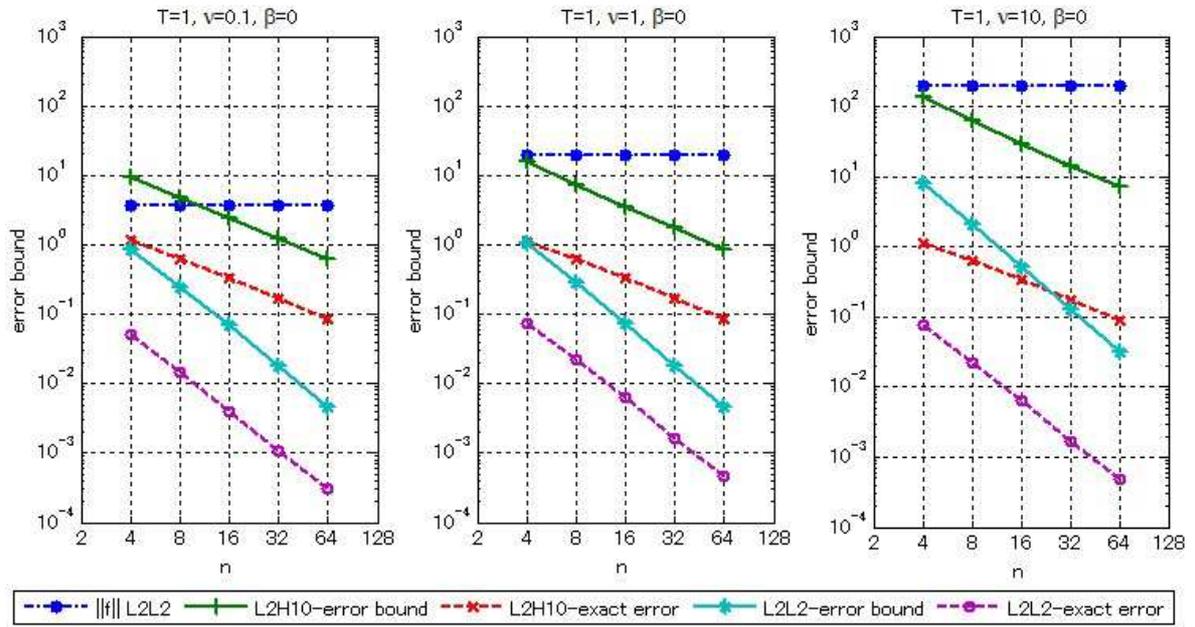}
\caption{~\text{L2H10 and L2L2 error estimates,~$m=n^2$,~~$\beta=0$,~~$T=1$,~~$\nu=0.1,~1,~10$,}}
\label{fig:1}
\end{center}
\end{figure}
\begin{figure}[htbp]
\begin{center}
\includegraphics[clip,width=0.95\linewidth]{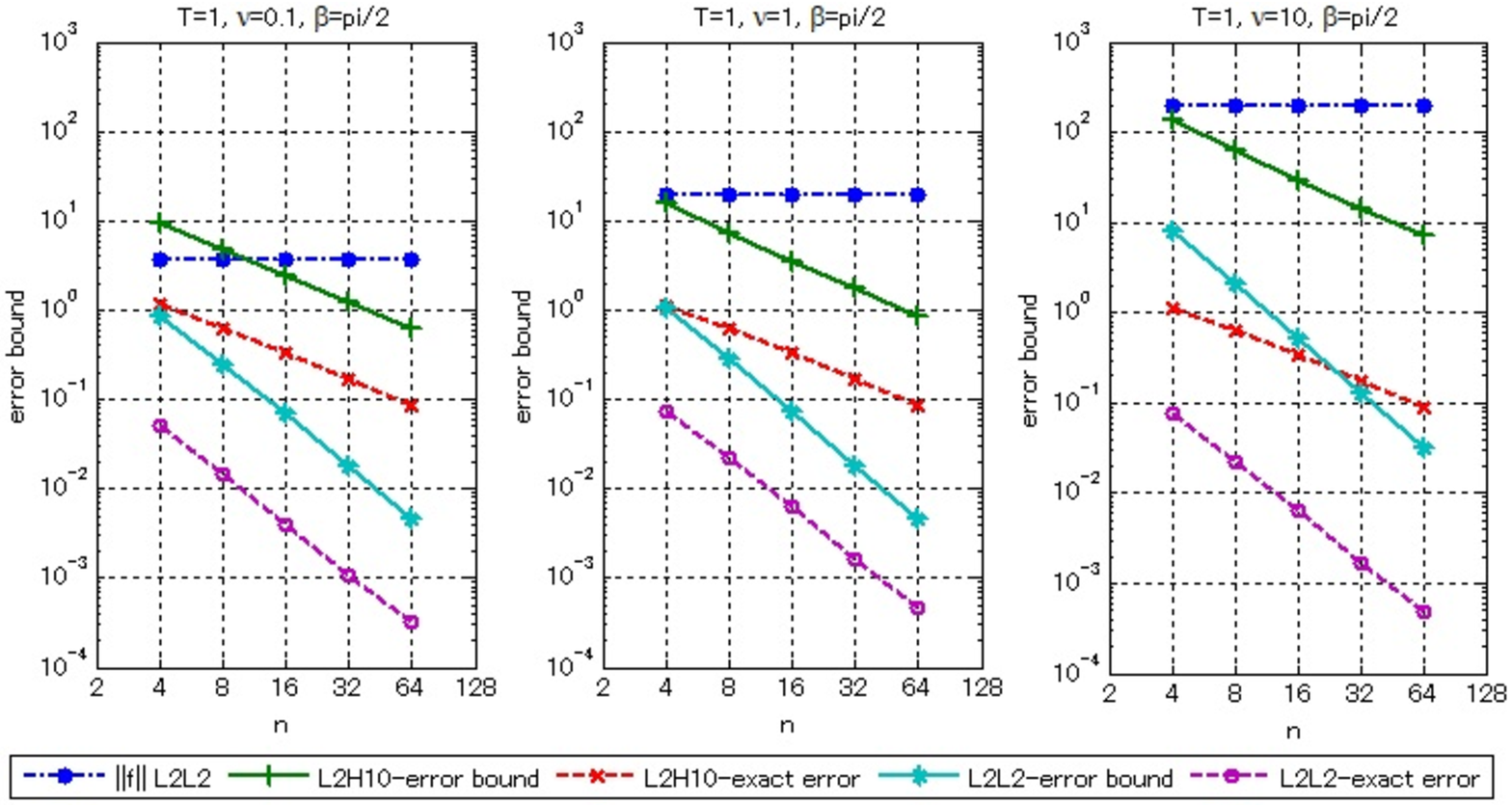}
\caption{~\text{L2H10 and L2L2 error estimates,~$m=n^2$,~~$\beta=0.5\pi$,~~$T=1$,~~$\nu=0.1,~1,~10$,}}
\label{fig:1}
\end{center}
\end{figure}




\begin{thebibliography}{00}
\bibitem{Barbu 1998}
V. Barbu, Partial differential equations and boundary value problems, \textit{Kluwer Academic Publishers, the Netherland}, 1998.

\bibitem{Bernardi 1982}
C. Bernardi, Numerical approximation of a periodic linear parabolic problem, SIAM  Journal on Numerial Analysis, 19 (1982), 1196-1207.

\bibitem{Figueras 2017}
J.-L. Figueras, M. Gameiro, J.-P. Lessard and R. de la Llave. A framework for the numerical computation and a posteriori verification of invariant objects of evolution equations. SIAM Journal on Applied Dynamical Systems 16(2), 1070-1088, 2017.

\bibitem{Hansbo 1991}
A. Hansbo, Error estimates for the numerical solution of a time-periodic linear parabolic problem, BIT 31, 664-685, 1991.

\bibitem{Nakao2011}
M.T. Nakao, T. Kinoshita and T. Kimura: 
On a posteriori estimates of inverse operators for linear parabolic initial-boundary value problems, 
 \textit{Computing}, 94, 151-162, 2012.
%
\bibitem{Nakao 2005}
M.T. Nakao, K. Hashimoto, Y. Watanabe: 
A numerical method to verify the invertibility of linear elliptic operators with applications to nonlinear problems, 
\textit{Computing }, {75}, 1--14, 2005.
%
\bibitem{Nakao SIAM}
 M. T. Nakao, T. Kimura, T. Kinoshita, Constructive a priori error estimates for a full discrete approximation of the heat equation, SIAM Journal on Numerial Analysis, 51 (2013), 1525-1541.


\bibitem{Pazy 1983}
A.Pazy, 
Semigroups of Linear Operators
and Applications to
Partial Differential Equations, \textit{Springer, New York}, 1983. 


\bibitem{INTLAB}
S.M. Rump, 
INTLAB - INTerval LABoratory. In Tibor Csendes, editor, Developments in Reliable Computing, pages 77-104. Kluwer Academic Publishers, Dordrecht, 1999. 
\bibitem{Schultz 1973}
M.H. Schultz,
Spline Analysis,
\textit{Prentice-Hall, N.J.}, 1973.
\bibitem{Zeidler 1990}
E. Zeidler, 
Nonlinear functional analysis and its applications II/B, 
\textit{Springer-Verlag, New York}, 1990.
\bibitem{Zgliczynski 2010}
P. Zgliczynski, 
{\em Rigorous Numerics for Dissipative PDEs III. An effective algorithm for rigorous integration of dissipative PDEs}, Topological Methods in Nonlinear Analysis, {36} (2010), pp.~197--262 
\end{thebibliography}



\end{document}